\documentclass[american]{article}
\usepackage{amssymb}
\usepackage{amsmath}
\newtheorem{all}{}

\newcommand{\mk}{\medskip}
\newcommand{\bk}{\bigskip}

\newcommand{\Nx}{{G_{[\,\overline{\Omega}\,]}}_X}
\newcommand{\No}{{G_{[\,\overline{\Omega}\,]}}_O}
\newcommand{\N}{G_{[\,\overline{\Omega}\,]}}

\newcommand{\B}{\hfill$\Box$}
\newcommand{\centre}{{\mathfrak z}(G_u)}
\newcommand{\punti}{,\,\ldots\, ,}

\newcommand{\puntini}{\,\ldots\,}

\newcommand{\car}{{\mathrm{char}}\,{\sf k}}

\newcommand{\ka}{{\sf k}}
\newcommand{\à}{\`a}

\newcommand{\é}{\'e}

\begin{document}
\title{Algebraic $(2,2)$-transformation groups\footnote{This paper contains the more significant part of the article with the same title that
will appear in the Volume 12 of Journal of Group Theory}}
\author{C.~Bartolone, A.~Di Bartolo, K.~Strambach}

\date{}

\large \sloppy \maketitle

\abstract{\footnotesize In this paper we determine all algebraic
transformation  groups $G$, defined over an algebraically closed
field $\sf k$, which operate transitively, but not primitively, on
a variety $\Omega$, provided the following conditions are
fulfilled. We ask that the (non-effective) action of $G$ on the
variety of blocks is sharply 2-transitive, as well as the action
on a block $\Delta$ of the normalizer $G_\Delta$. Also we
require sharp transitivity on pairs $(X,Y)$ of independent points
of $\Omega$, i.e. points contained in different blocks.}

\paragraph{} Although
classifications of imprimitive permutation groups appeared already
at beginning of the last century (see \cite{Martin}) and imprimitive actions
 play an
important role in geometry, the corresponding literature is
actually less well-developed than the one concerning primitive
groups. For finite groups some classification has been done
 (see for instance
\cite{Silvia}, \cite{Dixon} and \cite{Ronse}). In \cite{Silvia}
by using wreath products, the best-known construction principle to
get imprimitive groups, a classification  of finite
imprimitive groups, acting highly transitively on blocks and satisfying conditions
 very common in geometry, is achieved.

The present paper arises with the aim to obtain classifications
for
 infinite
 imprimitive groups belonging to well-studied categories. We start with an  imprimitive algebraic
 group
 $G$,
 over an algebraically closed field $\ka$, operating on an algebraic
variety $\Omega$ of positive dimension in such a way that the
induced actions on
 the set $\overline{\Omega}$ of
 blocks and on a block $\Delta$ are both
 sharply 2-transitive. Moreover we ask
  the group to act sharply transitively on pairs of points lying
 in different blocks. The latter condition,  frequently occurring in
geometry (see for instance \cite{Benz}),  avoids a too general
 context. For the classification we do not need the group actions
  be bi-regular morphisms but we just ask that the orbit maps be
 separable morphisms.
 It turns out that $G$ is
 the semidirect product of a 3-dimensional  unipotent connected
 group $G_u$ by a 1-dimensional connected torus $T$, both acting on the
points of an affine
 plane over $\ka$ with a full set of parallel lines as the blocks.

 There are
  two subgroups which play a fundamental role for the classification:
 the kernel $\N$ of the representation on $\overline{\Omega}$ (the
so-called {\it inertia
 subgroup}) and its stabilizer $\No$ of a fixed point $O$, which turns out
to be even the
 point-wise stabilizer of the block containing $O$. There
 exists a $G$-invariant transversal $L$ of $G$ with respect to $\No$
 which is essential for the classification. $L$ is a subgroup precisely if
 $G_u/\centre$ is commutative, in such a case $G_u/\centre$ is even a vector group. Fixing the
 structure
of $L$, the
 classification (see the main theorem) depends on four (not
 necessarily independent)
 integer parameters which distinguish the isomorphism class of $G$.
 But if the $\car$ is positive, then for suitable values of the integer parameters it happens that $L$ could
 be both a vector group and a non-commutative group.

We refer to \cite{Rosenlicht} for well-known results about
non-affine algebraic groups and to \cite{Humphreys} about affine
algebraic groups.

\paragraph{$\bf \S1.$} Throughout the paper $G$ will denote an
algebraic group defined over an algebraically closed field $\sf
k$, operating effectively on the points of a  variety $\Omega$ of
positive dimension. We assume that the orbit maps $g\mapsto g(X)$
are separable morphisms $G\to\Omega$ and $G$ acts transitively
with a nontrivial system of imprimitivity $\overline{\Omega}$.
Moreover, putting
\begin{description}
    \item[-] the {\it normalizer} $G_\Delta:=\left\{g\in
G:g(\Delta)=\Delta\right\}$ of $\Delta\in\overline{\Omega}$,
    \item[-] the {\it centralizer} $G_{[\,\Delta\,]}:=\left\{g\in
G_\Delta:g(X)=X\;\forall
    X\in\Delta\right\}$ of $\Delta\in\overline{\Omega}$,
    \item[-] the   {\it inertia} subgroup
$G_{[\,\overline{\Omega}\,]}:=\left\{g\in
    G:g(\Delta)=\Delta\;\forall\Delta\in\overline{\Omega}\right\}$,
\end{description}
we require the following transitivities:
\begin{itemize}
    \item[1.] $G_\Delta/G_{[\,\Delta\,]}$ acts sharply 2-transitively on
    $\Delta$,
    \item[2.] $G/G_{[\,\overline{\Omega}\,]}$ acts sharply 2-transitively on
    $\overline{\Omega}$,
    \item[3.] $G$ acts sharply transitively on
    $\Lambda:=\left\{(X,Y)\in\Omega^2:\Delta_X\ne\Delta_Y\right\}$, where $\Delta_Z\in\overline{\Omega}$ denotes the block containing
$Z\in\Omega$.
\end{itemize}
 We call such a triple ${\sf
G}=(G,\Omega,\overline{\Omega})$ a $(2,2)$-{\it imprimitive
algebraic group}. Since the stabilizer of a point is not trivial,
 conditions 3 and 1 guarantiy that the centre of $G$ consists just of the
identity. Hence the algebraic group $G$ must be affine.

\begin{all}{\bf $\!\!$. Proposition:}\label{2}
\begin{enumerate}
    \item[$i)$] Every block $\Delta\in\overline{\Omega}$ is
    closed and $G_\Delta=\N G_X$ for any $X\in\Delta$;
    \item[$ii)$] the inertia subgroup $\N$ is closed.
\end{enumerate}
\end{all}

\begin{all}{$\!\!${\bf .}$\;$\it Remark}\label{1}{
\rm : As  orbit maps are separable morphisms $G\to\Omega$, by the
universal mapping property we may identify $\Omega$ with the
homogeneous space $G/G_O$ for a fixed stabilizer $G_O=\big\{g\in
G:g(O)=O\big\}$, $O\in\Omega$. As well as, in view of Proposition
\ref{2}, we may identify $\overline{\Omega}$ with the homogeneous
space $G/G_\Delta$.}
\end{all}

\begin{all}{\bf $\!\!$. Proposition:}\label{4}
For all $X\in\Omega$ the centralizer
$G_{{[\,\overline{\Omega}\,]}_X}=\big\{g\in\N:g(X)=X\big\}$ is
contained in $G_{[\,\Delta_X]}$ and
$G_{[\,\overline{\Omega}\,]}=G_{{[\,\overline{\Omega}\,]}_X}\!\times
G_{{[\,\overline{\Omega}\,]}_Y}$ for any $(X,Y)\in\Lambda$.
\end{all}

\begin{all}{\bf $\!\!$. Proposition:}\label{Delta}
\begin{description}
    \item[$a)$] $\overline{\Omega}$ contains infinitely many
    blocks and
every block contains infinitely many points;
     \item[$b)$] $G_O$ is the semidirect product of the $1$-dimensional
    connected unipotent subgroup $\No$ by a $1$-dimensional connected torus $T$;
    \item[$c)$]  $G/\N$ is a $2$-dimensional Frobenius
algebraic group with complement $\simeq T$; \item[$d)$] For all
$\Delta\in\overline{\Omega}$, $G_\Delta/G_{[\,\Delta\,]}$ is a
$2$-dimensional Frobenius algebraic group whose $1$-dimensional
kernel is isomorphic to $\Nx$ for any $X\in\Omega\setminus\Delta$.
\end{description}
  \end{all}

\begin{all}{\bf $\!\!$. Proposition:}\label{5}
$G$ is a solvable connected affine group of dimension $4$ and $G$
is the semidirect product of its unipotent radical $G_u$ by the
torus $T$. Moreover the centre ${\mathfrak z}(G_u)$ of $G_u$ is
contained in $\N$ and for any $X\in\Omega$ we have $\N={\mathfrak
z}(G_u)\times\Nx$.
\end{all}

\begin{all}{$\!\!${\bf .}$\;$\it Remark}\label{8}{
\rm : If we denote  by $g_u$ and $g_s$ the images of $g\in G$
under the projections $G_u\times T\to G_u$ and $G_u\times T\to T$,
respectively, the mapping $\pi:G\to G_u/\No$ with $\pi(g)=g_u\No$
turns out to be a separable morphism of algebraic varieties. The
fibres of $\pi$ are precisely the cosets $gG_O$, so $gG_O\mapsto
g_u\No$ yields an isomorphism $G/G_O\to G_u/\No$. So we may take
the homogeneous space $G_u/\No$ as $\Omega$ and
$$
\Big(g,h\No\Big)\mapsto ghg_s^{-1}\No\quad (g\in G,h\in G_u)
$$
as the action of $G$ on $\Omega$ since
$(g_1g_2)_u=(g_1)_u(g_1)_s(g_2)_u(g_1)_s^{-1}$. In particular
$\Omega\simeq G_u/\No$ is a 2-dimensional (irreducible affine)
variety with

\centerline{$\overline{\Omega}=\bigcup_{g\in
G_u}\!\!\Delta_{g(O)}\simeq\bigcup_{g\in G_u}g\centre\No.$}}
\end{all}

\paragraph{$\bf \S2.$} Let $G=U\rtimes T$ be a semidirect product of an $n$-dimensional connected unipotent group $U$ by a 1-dimensional
connected torus $T$. According to Serre \cite{Serre}, p. 172, the
group $U$ has a representation on the affine space  $\ka^n$ in
such a way the subspaces $$U_i=\big\{(x_1\punti
x_n)\in\ka^n:x_{i+1}=\puntini =x_n=0)\big\}$$ are normal subgroups
of $G$, the product is given by $(x_1\punti x_n)(y_1\punti y_n)=$
$$
\big(x_1+y_1+\psi_1(x_2\punti x_n,y_2\punti y_n)\punti
x_{n-1}+y_{n-1}+\psi_{n-1}(x_n,y_n),\, x_n+y_n\!\big),
$$ for suitable polynomials $\psi_j\in\ka[x_j\punti x_n,y_{j+1}\punti
y_n]$, and the automorphism of $U$ induced by an element $\tau\in
T$ maps $(x_1\punti x_n)$ to
$$
\big(a_\tau^{e_1}x_1+\varphi^{(\tau)}_{1}(x_2\punti x_{n}\!)\punti
a_\tau^{e_{n-1}}x_{n-1}+\varphi_{n-1}^{(\tau)}(x_n),\,
a_\tau^{e_n}x_n\!\big)$$ with  $a_\tau\in{\sf k}^\ast$, an element
depending bi-regularly on $\tau$, the map $\varphi_j^{(\tau)}$ a
morphism $U_n/U_j\to U_j/U_{j-1}$ and $e_j$ a fixed integer.

\begin{all}{$\!\!$\bf .\;Lemma:}\label{a.2} Let $n\geq 2$. Then for any
$\tau\in T$ the morphism $\varphi^{(\tau)}_{n-1}$ yields a group
homomorphism $U_n/U_{n-1}\to U_{n-1}/U_{n-2}$. Moreover we may
take as $\psi_{n-1}$
$$
  \begin{array}{ll}
  a)\;\hbox{the zero polynomial}, & \hbox{if $U_n/U_{n-2}$ is a vector
group,} \\
   b)\;\sum_{i=1}^{p-1}\frac
1p\left(^{p\,}_{\,i}\right)x_n^{i{p^r}}y_n^{(p-i)p^{r}}, &
\hbox{if
$U_n/U_{n-2}$ is an Abelian group of exponent $p^2$,}\\
    c)\;x_n^{p^r}y_n^{p^{s}}, & \hbox{if $U_n/U_{n-2}$ is not commutative,}
  \end{array}
$$
where, in the cases $b)$ and $c)$, $p=\car >0$, $r$ and $s$ are
nonnegative integers such that $r<s$ and $e_{n-1}=e_n\!\deg
(\psi_{n-1})$.
\end{all}

\begin{all}{$\!\!${\bf .}$\;$\it Remark}\label{dieci}{
\rm : It follows from \cite{Birula} that the action of a
1-dimensional torus on a 2-dimensional connected unipotent group
$U$ may be given by diagonal $(2\times 2)-$matrices with entries
in $\ka$. The following lemma, which generalizes both the lemma on
p. 109 in \cite{Rosenlicht2} and Corollary 2.9 in \cite{Falcone},
shows that this can be done without destroying the group structure
of $U$.} \end{all}

\begin{all}{$\!\!$\bf .\;Lemma:}\label{undici}
Let $\varphi_2^{(\tau)}=\puntini=\varphi_{n-1}^{(\tau)}=0$ and
assume $\varphi_{1}^{(\tau)}$ is  a group homomorphism
$U_n/U_{n-1}\to U_1$. Then there exists a bi-regular section
$\sigma:U_n/U_{n-1}\to U_n$ such that $\sigma
\big(x_nU_{n-1}\big)=\big(f(x_n),0\punti 0,x_n\big)$ with
$\delta^1\!(f)=0$ and  $\sigma\big(U_n/U_{n-1}\big)$ invariant
under $T$.
\end{all}

Set $M:=\big\{(0\punti 0,x_n):x_n\in\ka\big\}$ and let
$\upsilon=(0\punti 0,u)\in M$. We have

\begin{all}{\bf $\!\!$. $\!\!$Lemma:}\label{10}  Let $n\geq 3$.
Assume the centralizer ${\mathfrak C}_{U_{n-1}}(\upsilon)$ of
$\upsilon$ in $U_{n-1}$ satisfies the condition ${\mathfrak
C}_{U_{n-1}}(\upsilon)=U_{n-2}\mod U_{n-3}$ for all $\upsilon\in
M$. Then the automorphism $\rho_\upsilon$ of $U_{n-1}/U_{n-3}$
induced by conjugation by $\upsilon$ maps
$$(\puntini,x_{n-2},
x_{n-1},0)U_{n-3}\longmapsto\big(\puntini , x_{n-2}+
u^hx_{n-1}^{k},x_{n-1},0\big)U_{n-3}$$ with  $h$ and $k$ p-powers
if $\car =p>0$, $h=k=1$ otherwise.
\end{all}

In the remaining part of the paper we ask the torus $T$ to act
sharply transitively on $U_n/U_{n-1}$. This
means\begin{equation}\label{e1}
    e_n=\left\{%
\begin{array}{ll}
    1, & \hbox{if $\mathrm{char}\,{\sf k}=0$;} \\
    \hbox{a $p$-power,} & \hbox{if $\mathrm{char}\,{\sf k}=p>0$.} \\
\end{array}%
\right.
\end{equation}

\paragraph{$\bf \S3.$} Now we go back to the  the  $(2,2)$-{imprimitive
algebraic group} ${\sf G}=(G,\Omega,\overline{\Omega})$. This
section is devoted to the case where the 2-dimensional factor
group $G_u/{\mathfrak z}(G_u)$ is {\it commutative}.

\begin{all}{\bf $\!\!$. Proposition:}\label{13}
$G_u/{\mathfrak z}(G_u)$ is a vector group.
\end{all}

\begin{all}{\bf $\!\!$. Proposition:}\label{14P} There exists a
$T\!$-invariant normal subgroup $L$ of $G_u$ containing the centre
${\mathfrak z}(G_u)$ and $G_u=L\rtimes\No$.
\end{all}

According to the notation of \S2 we may take $U_1=\centre$,
$U_2=\N$, $U_3=G_u$. In addition we may choose
$$\No=\big\{(0,x_2,0): x_2\in {\sf k}\big\},$$ the subgroup $\No$ being $T$-invariant. Observing that the
normal subgroup $L$ of $G$ is not contained in $U_2$, we may also
put
$$
L=\big\{(x_1,0,x_3): x_1, x_3\in {\sf k}\big\}.
$$
Thus  the product $(x_1,0,x_3)(y_1,0,y_3)$ of two elements of $L$
is  given by
$$\big(x_1+y_1+\beta(x_3, y_3), 0, x_3+y_3\big)$$ and by Lemma \ref{a.2}
we may take {
\begin{equation}\label{IV}
\beta(x_3, y_3)=\left\{%
\!\!\!\begin{array}{ll}
    0, & \!\!\!\!\hbox{if $L$ is a vector group,} \\
    \sum_{i=1}^{p-1}\frac 1
p\big(^{p\,}_{\,i}\big)x_3^{i{p^r}}y_3^{(p-i)p^{r}}\!, &
\!\!\!\!\hbox{if
$L$ is commutative of exponent $p^2$,} \\
   x_3^{p^r}y_3^{p^{s}}, & \!\!\!\!\hbox{if $L$ is not commutative,} \\
\end{array}%
\right. \end{equation} }  $\!\!$for some nonnegative integers
$r,s$ with $r<s$. Besides an element ${\upsilon}=(0,0,u)\in L$
moves the block $\Delta_O$ to a different block
$\Delta_{\upsilon(O)}$ (Remark \ref{8}), so $\upsilon$ centralizes
no element in $\No$, the intersection $\No\cap{\N}_{\ell(O)}$
being trivial. Then Lemma \ref{10} applies and, up to the
isomorphism $(x_1,x_2,x_3)\mapsto(x_1,c^{\frac 1{h_2}}x_2,x_3)$,
we may claim

\begin{all}{\bf $\bf\!\!.$ Proposition:}\label{15}
The product $(x_1,x_2,x_3)(y_1,y_2,y_3)$ in $G_u$ may be defined
through
$$\left(x_1+y_1+y_2^{h_2}x_3^{h_3}+\beta(x_3,
y_3),\, x_2+y_2,\, x_3+y_3\right),$$ where $\beta$ is given by
{\rm (\ref{IV})} and each exponent $h_i$ is a $p$-power in case
$\mathrm{char}\,{\sf k}=p>0$, $h_i=1$ otherwise.\B
\end{all}

As we observed in Remark \ref{dieci}, there is no loss of
generality if we assume the action of the torus $T$ on the affine
plane $L$ given by diagonal $(2\times 2)-$matrices. But $\No$
occurs as a further $T$-invariant subgroup of dimension 1, so the
diagonal action of each $\tau\in T$ extends to the whole group
$G_u$ via
\begin{equation}\label{tau}
(x_1,x_2, x_3)\mapsto\big(a_\tau^{e_1} x_1,\,a_\tau^{e_2}
x_2,\,a_\tau^{e_3} x_3\big).
\end{equation}
The value of the exponent $e_3$ was given by (\ref{e1}), whereas
the possible relationship occurring between $e_1$ and $e_3$ was
stated in Lemma \ref{a.2}. Now by imposing that $\tau$ is a group
homomorphism we find
\begin{equation}\label{e3}
    e_1=e_2h_2+e_3h_3
\end{equation}
with $h_i$ arising from the product of $G_u$ given in Proposition
\ref{15}.

\paragraph{$\bf \S4.$} Assume now the factor
group $G_u/\centre$ to be {\it not commutative}. This requires
$\car =p>0$ and we are going to see that even $p>2$ holds.

Referring to the notation of \S2 we may take again $U_3=G_u$,
$U_2=\N$, $U_1=\centre$ and
$$
\No=\big\{(0,x_2,0): x_2\in {\sf k}\big\}.
$$
Also, by Lemma \ref{a.2},
$$\psi_2:(x_3,y_3)\mapsto x_3^{p^m}y_3^{p^{n}},$$ for some integer
$p$-powers $p^m$ and $p^n$ such that $m<n$. Furthermore,
 looking at Remark \ref{8}, we see that an element
$\upsilon=(0,0,x_3)$ moves the block $\Delta_O$ to a different
block $\Delta_{\upsilon(O)}$. So $\upsilon$ does not centralize
any element of $\No$ because the intersection
$\No\cap{\N}_{\upsilon(O)}$ is assumed to be trivial. So Lemma
\ref{10} applies and, up to an isomorphism, we may assume that the
 automorphism induced on $\N$ by an element $(0,0,x_3)$ maps
$$(y_1,y_2,0)\mapsto\big(y_1+y_2^{h_2}x_3^{h_3},y_2,0)$$
     for suitable integer $p$-powers $h_i=p^{l_i}$, $i=2,3$. If we represent
     $G_u$ as a noncentral extension of the vector group $\N$ by $G_u/\N$ using the cross section
$(x_1,x_2,x_3)\N\mapsto(0,0,x_3)$, the product
$(x_1,x_2,x_3)(y_1,y_2,y_3)$ of two elements in $G_u$ can also be
given by
$$
(x_1+y_1+y_2^{h_2}x_3^{h_3}+\beta(x_3,y_3),
x_2+y_2+x_3^{p^m}y_3^{p^n},x_3+y_3)
$$ with $\beta$ in $\ka[x_3,y_3]$ such that
$\beta(0,y_3)=\beta(x_3,0)=0$ and $G_u$ is determined by taking
$$
\psi_1(x_1,x_2,y_1,y_2)=y_2^{h_2}x_3^{h_3}+\beta(x_3,y_3).
$$
Now associative law forces the polynomial
$$\delta^2(\beta)(z_1,z_2,z_3)=\beta(z_1, z_2)+\beta(z_1+z_2,
z_3)-\beta(z_2, z_3)-\beta(z_1, z_2+z_3)$$ to be
\begin{equation}\label{I}\delta^2(\beta)(z_1,z_2,z_3)=
z_1^{p^{l_3}}z_2^{p^{l_2+m}}z_3^{p^{l_2+n}}\end{equation} and we
can state

\begin{all}{\bf $\bf\!\!.$ Proposition:}\label{18}
A necessary and sufficient condition in order that $G_u$ can be
constructed as an extension of $\centre$ by a non-commutative
connected unipotent group is that there exists a polynomial
$\beta\in\ka[x_3,y_3]$ satisfying {\rm (\ref{I})} with
$\beta(0,y_3)=\beta(x_3,0)=0$. In such a case we may take
$\psi_1(x_2,x_3,y_2,y_3)=y_2^{p^{l_2}}x_3^{p^{l_3}}+\beta(x_3,
y_3)$.\B
\end{all}

The crucial question  now is under what conditions such a
polynomial $\beta$ there exists. Using a universal property of the
operator $\delta^2$ we have
$$\sum_{\pi\in{\sf
S}_n}\!\!\!{\mathrm{sign}(\pi)}\,\delta^2(\beta)(z_{\pi(1)},
z_{\pi(2)}, z_{\pi(3)})=0$$
\noindent and this, in view of
(\ref{I}), is equivalent to
\begin{equation}\label{l}
    \text{$l_3-l_2=m,\;$ or $\;l_3-l_2=n$.}
\end{equation}
Assume now $\mathrm{char}\,{\sf k}=2$ and ${l_2+m}>0$ and denote
by $\beta_j$ the homogeneous component of $\beta$ of degree $j$.
As in our case the operator $\delta^2$ is additive, (\ref{I}) says
that $\delta^2(\beta)=\delta^2(\beta_k)$, where
$k=2^{l_2}(2^q+2^m+2^n)$ with either $q=m$, or $q=n$ according as
whether $l_3-l_2=m$, or $l_3-l_2=n$. Let

\centerline{$\beta_{k}(y_1,y_2)=\sum_{i=0}^{k}a_iy_1^{k-i}y_2^i.$}

\noindent Then (\ref{I}) becomes
$$
\begin{array}{l}
\sum_{i=0}^{k}\!a_i\!\left(z_1^{k-i}z_2^i+(z_1+z_2)^{k-i}z_3^i+z_2^{k-i}z_3^i
   +z_1^{k-i}(z_2+z_3)^i\right)=z_1^{2^{l_2+q}}
z_2^{2^{l_2+m}}z_3^{{2}^{l_2+n}}. \\
\end{array}
$$
\noindent Deriving this identity with respect to $z_1$ and
evaluating at $(0,y_1,y_2)$ we obtain
\begin{equation}\label{aaaa}
a_{k-1}y_1^{k-1}+\frac\partial{\partial
y_1}\beta_{k}(y_1,y_2)+a_{k-1}(y_1+y_2)^{k-1}=0,
\end{equation}
whereas deriving with respect to $z_3$ and evaluating at
$(y_1,y_2,0)$ we get
\begin{equation}\label{aaaaa}
a_1(y_1+y_2)^{k-1}+a_1y_2^{k-1}+\frac\partial{\partial
y_2}\beta_{k}(y_1,y_2)=0.
\end{equation}

\noindent As $\car=2$, $\frac\partial{\partial
y_1}\beta_{k}(y_1,y_2)$ and $\frac\partial{\partial
y_2}\beta_{k}(y_1,y_2)$ are polynomials in $y_1^2$ and $y_2^2$,
respectively, the identities (\ref{aaaa}) and (\ref{aaaaa}) force
$a_{k-1}=a_1=0$, hence
$$\frac\partial{\partial y_1}\beta(y_1,y_2)=\frac\partial{\partial
y_2}\beta(y_1,y_2)=0$$ and this yields $a_i=0$ for all odd $i$.
Thus we may do the substitution
$(z_1,z_2,z_3)\mapsto(z_1^2,z_2^2,z_3^2)$, hence
$(z_1,z_2,z_3)\mapsto(z_1^{2^{l_2+m}},z_2^{2^{l_2+m}},z_3^{2^{l_2+m}})$
by iterating the process. So (\ref{I}) turns into
\begin{equation}\label{II}
\delta^2(\gamma)(z_1,z_2,z_3)=z_1^{2^{q-m}}z_2z_3^{2^{n-m}}
\end{equation}
with $\gamma\big(y_1^{2^{l_2+m}},y_2^{2^{l_2+m}}\big)=\beta
(y_1,y_2)$. Let $\gamma_t$ be the homogeneous component of degree
$t:=1+2^{n-m}+2^{q-m}$ of $\gamma$ and let
$\gamma_t(y_1,y_2)=\sum_{i=0}^{t}b_iy_1^{t-i}y_2^i$. Then
(\ref{II}) says that $\delta^2(\gamma)=\delta^2(\gamma_t)$, hence
$$\sum_{i=0}^{t}\!b_i\!\left(z_1^{t-i}z_2^i+(z_1+z_2)^{t-i}z_3^i+z_2^{t-i}z_3^i
   +z_1^{t-i}(z_2+z_3)^i\right)=z_1^{2^{q-m}}z_2z_3^{2^{n-m}}\!.$$

\noindent Likewise above we obtain
$$
\left\{%
\begin{array}{ll}
    b_{t-1}y_1^{t-1}+\frac\partial{\partial
y_1}\gamma_t(y_1,y_2)+b_{t-1}(y_1+y_2)^{t-1}=(1-\epsilon)y_1y_2^{2^{n-m}},\\
\\
    b_1(y_1+y_2)^{t-1}+b_1y_2^{t-1}+\frac\partial{\partial
y_2}\gamma_t(y_1,y_2)=0,
\end{array}%
\right.
$$ where either $\epsilon=0$, or $\epsilon =1$ according as
whether $q=m$, or $q=n$. So  Euler's identity says that
$t\gamma_t(y_1,y_2)$ is the polynomial{\footnotesize
$$
b_{t-1}y_1\!\!\left(\!y_1^{t-\!1}\!\!+\!
(y_1+y_2)^{t-\!1}\!\!\!+(1-\epsilon)y_1y_2^{2^{n-m}}\right)\!+b_1y_2\!\!\left(\!(y_1\!+y_2)^{t-\!1}\!\!+y_2^{t-\!1}\!\right)
$$}
$\!\!$or the polynomial {\footnotesize
$$
\begin{array}{l}
b_{t-1}\left(y_1^{1+2^{n-m}}y_2^{2^{q-m}}+
y_1^{1+2^{q-m}}y_2^{2^{n-m}}+y_1y_2^{2^{n-m}+2^{q-m}}+(1-\epsilon)\,y_1^2y_2^{2^{n-m}}\right)+\qquad\quad\\
\\
\hfill + b_1\left(y_1^{2^{n-m}+2^{q-m}}y_2+ \
y_1^{2^{n-m}}y_2^{1+2^{q-m}}+y_1^{2^{q-m}}y_2^{1+2^{n-m}}\right).
\end{array}
$$}
Let $q=m$. Then we have the polynomial identity
$$
\big(b_{t-1}+b_1\big)\left(y_1^{1+2^{n-m}}y_2+y_1y_2^{1+2^{n-m}}\right)+b_{1}y_1^{2^{n-m}}y_2^2=0
$$
which asks $b_{t-1}=b_1=0$ and, consequently,
$\frac\partial{\partial y_1}\gamma_t(y_1,y_2)=y_1y_2^{2^{n-m}}$, a
contradiction. Let $q=n$. Then
$$
\gamma_t(y_1,y_2)=b_{t-1}y_1y_2^{2^{n-m+1}}+b_1y_1^{2^{n-m+1}}\!\!y_2
$$
and $\delta^2(\gamma_t)(x_1,x_2,x_3)=0$. This contradicts
(\ref{II}) and $p\ne 2$ follows.

Actually, if $p\ne 2$ the polynomials
\begin{equation}\label{25.5.07}\beta(x_3,y_3)=\begin{cases}
\frac 12 x_3^{2p^{l_3}}y_3^{p^{l_2+n}}& \text{if}\,\,
l_3-l_2=m;\\
x_3^{p^{l_3}+p^{l_2+m}}y_3^{p^{l_3}}\!\!+\frac 12
x_3^{p^{l_2+m}}y_3^{2p^{l_3}}& \text{if}\,\, l_3-l_2=n.
\end{cases} \end{equation} satisfy the
conditions required in Proposition \ref{18}. Any other polynomial
satisfying the conditions of Proposition \ref{18} differs from
(\ref{25.5.07}) for a co-cycle $\kappa(x_3,y_3)$ for a central
extension of ${\sf k}_+$ by ${\sf k}_+$ that we are going to show
it is a co-boundary.

It follows from Remark \ref{dieci} that we may assume any element
$\tau\in T$  acts on $G_u$ via
$$(x_1,x_2,x_3)\mapsto \big(a_\tau^{e_1}x_1+
\varphi_1^{(\tau)}(x_3),a_\tau^{e_2}x_2, a_\tau^{e_3}x_3\big),$$
with the morphism $\varphi_1^{(\tau)}$ depending only on $x_3$
because $\No$ is $T$-invariant. By imposing that $\tau$ operates
as a group homomorphism we obtain first
\begin{equation}\label{e3bis}
\hbox{$e_2=e_3(p^m+p^n)$ and $e_1=e_2 h_2+e_3h_3
=e_3\big(p^{l_3}+p^{l_2+m}+p^{l_2+n}\big)$},
\end{equation}
but also
$$a_\tau^{e_1}
\beta(x_3, y_3)-\beta\big(a_\tau^{e_3} x_3, a_\tau ^{e_3}
y_3\big)+a_\tau^{e_1} \kappa(x_3, y_3)- \kappa\big(a_\tau^{e_3}
x_3, a_\tau ^{e_3}
y_3\big)=\delta^1\!\big(\varphi_1^{(\tau)}\big)(x_3,y_3),$$ or
\begin{equation}\label{kappa}
    a_\tau^{e_1}
\kappa(x_3, y_3)-\kappa\big(a_\tau^{e_3} x_3, a_\tau ^{e_3}
y_3\big)=\delta^1\!\big(\varphi_1^{(\tau)}\big)(x_3,y_3),
\end{equation}
because $e_1=e_3\!\deg\beta$ in view of (\ref{e3bis}). Since $e_3$
is a $p$-power and $p>2$, the integer $e_1$ can be, by
(\ref{e3bis}), neither a $p$-power, nor the sum of two $p$-powers.
Thus Theorem 4.6 in \cite{DMG} guaranties that $\kappa$ is a
co-boundary, i.e. $\kappa=\delta^1\!(g)$ for some polynomial
$g\in\ka[\mathrm{T}]$, that may be eliminated using the
substitution $x_1\mapsto x_1-g(x_3)$. Such a replacement yields
$\delta^1\!\big(\varphi_1^{(\tau)}\big)(x_3,y_3)=0$, i.e.
$\varphi_1^{(\tau)}$ is additive, and we may assume the action of
$T$ given by diagonal matrices, as Lemma \ref{undici} claims.

\paragraph{$\bf \S5.$} Now we  collect all information  achieved
in the previous sections and classify  $\sf G$ according to the
structure of the transversal $L$. With the aid of Remark \ref{8}
we can state:

\begin{all}{\bf $\bf\!\!.$ Main Theorem:}\label{twenty}
Every $(2,2)$-imprimitive algebraic group ${\sf
G}=(G,\Omega,\overline{\Omega})$ can be constructed on the affine
variety ${\sf k}^3\times{\sf k}^\ast$ as follows:
\begin{itemize}
    \item define the unipotent radical $G_u$ on the affine space
    $\ka^3$ through the
    product$$(x_1,x_2,x_3)(y_1,y_2,y_3)=\big(x_1+y_1+\psi_1(x_3,y_2,y_3),x_2+y_2+\psi_2(x_3,y_3),x_3+y_3\big),
    $$ where either
\begin{itemize}
    \item[] $\psi_2(x_3,y_3)=0$ and $\psi_1(x_3,y_2,y_3)=y_2^{h_2}x_3^{h_3}+\beta(x_3,y_3)$ with
each $h_i$ an integer $p$-power $p^{l_i}$ in case $\car=p>0$,
$h_i=1$ otherwise, and $\beta(x_3,y_3)$  one of the polynomials
\begin{itemize}
\item[$-$]  $0$; \item[$-$] $\sum_{i=1}^{p-1}\frac 1
p\big(^{p\,}_{\,i}\big)x_3^{i{p^r}}y_3^{(p-i)p^{r}}$; \item[$-$]
$x_3^{p^r}y_3^{p^s}$;
\end{itemize}
for suitable nonnegative integers $r,s$ such that $r<s$,
\end{itemize}
or
    \begin{itemize}
\item[] $\psi_2(x_3,y_3)=x_3^{p^m}y_3^{p^n}$, with $p=\car>2$
    and $m,n$ non-negative integers such that $m<n$, and
    $\psi_1(x_3,y_2,y_3)$  as above with
$$
\beta(x_3,y_3)=\left\{%
\begin{array}{ll}
\frac 12 x_3^{2p^{l_3}}y_3^{p^{l_2+n}}, &\hbox{if
$l_3-l_2=m$};\\
x_3^{p^{l_3}+p^{l_2+m}}y_3^{p^{l_3}}+\frac 12
x_3^{p^{l_2+m}}y_3^{2p^{l_3}}, &\hbox{if $l_3-l_2=n$};
\end{array}%
\right.$$
\end{itemize}
    \item leave $a\in\ka^\ast$ operate on $\ka^3$ via
    $$(u_1,u_2,u_3)\mapsto
    \big(a^{e_1}u_1,a^{e_2}u_2,a^{e_3}u_3\big)$$
    where \begin{itemize}
          \item $e_1=e_2h_2+e_3h_3$, but also
$e_1=e_3\!\deg\beta$ if $\beta$ is not the zero polynomial;
\item
$e_2=e_3\frac{\deg\beta-h_3}{h_2}$ if $\beta$ is not the zero
polynomial;
  \item
$e_3$ is a positive integer $p$-power in case $\car =p>0$,
    $e_3=1$ otherwise;
 \end{itemize}
    \item identify
$\Omega$ with the affine plane $\ka^2$ with the parallel lines
$y=k$ giving the set $\overline{\Omega}$ of blocks. Then a
transformation $(u_1,u_2,u_3,a)\in G$ moves the point
$(x,y)\in\Omega$ to the point
$$\big(u_1+a^{e_2h_2+e_3h_3}x+\psi_1(u_3,0,a^{e_3}y),\,
u_3+a^{e_3}y)\big).$$ \hfill$\Box$
\end{itemize}
\end{all}

The canonical representation of $\sf G$ given through the main
theorem  depends on  the polynomial $\beta$ as well as on the
integer parameters $e_2,e_3,h_2,h_3$, though $h_2$ and $h_3$ could
already be determined by $\beta$, $e_2$ and $e_3$. Labelling $\sf
G$ as ${\sf G}^{(e_2\!,\,e_3,\,h_2\!,\,h_3)}_\beta$, we ask
whether an  isomorphism
$$\Phi:{\sf G}^{(e_2\!,\,e_3,\,h_2\!,\,h_3)}_\beta\rightarrow {\sf
G}^{\big(e'_2\!,\,e'_3,h'_2\!,\,h'_3\big)}_{\beta'}$$ between two
$(2,2)$-imprimitive algebraic groups with different parameters
exists. Of course we may assume the same sets of points and blocks
for both  groups, so $\Phi$ is a pair $(\Phi_1,\Phi_2)$ with
$\Phi_1$ a group isomorphism
$G^{(e_2\!,\,e_3,\,h_2\!,\,h_3)}_\beta\longrightarrow
G^{\big(e'_2\!,\,e'_3,h'_2\!,\,h'_3\big)}_{\beta'}$ and
$\Phi_2:\ka^2\to\ka^2$ a bijective morphism of the affine plane
$\ka^2$ transforming horizontal lines into horizontal lines such
that
$$\Phi_2\big(g(P)\big)=\Phi_1(g)\big(\Phi_2(P)\big)\quad \big(g\in
G^{(e_2\!,\,e_3,\,h_2\!,\,h_3)}_\beta, P\in\ka^2\big).$$ As $G_u$
is transitive on $\Omega$, up to inner automorphisms we may assume
that $\Phi_2$ leaves the point $O=(0,0)$ of $\Omega$ fixed, hence
the line $y=0$ stable. Then the stabilizer of $O$, as well as the
normalizer and centralizers of $\Delta_O$ correspond; in
particular
\begin{equation}\label{cinque}
    \begin{array}{rlccccc}
    &\Phi_1\big((0,u_2,0)\big)=(0,b_2u_2,0)&(b_2\in\ka^\ast),&& \\
&\Phi_1\big((u_1,0,0)\big)=(b_1u_1,0,0)&(b_1\in\ka^\ast),&\qquad\qquad\qquad&\qquad&\quad\;&
    \end{array}
\end{equation}
and, moreover,
\begin{equation}\label{eq2}
\begin{array}{rlcccc}
    &\Phi_1\big((0,0, u_3)\big)=\big(f_1(u_3),f_2(u_3),
b_3u_3)\big)&(b_3\in\ka^\ast),
&\qquad\qquad&\qquad& \\
&\Phi_2\big((x,y)\big)= \big(b_1x+f_1(y),\,b_3y\big),&&
    \end{array}
\end{equation}
for suitable polynomials $f_j\in\ka[\mathrm{T}]$ such that
\begin{equation}\label{f_3}
\begin{array}{rlr}
  &
  \delta^1(f_2)(x_3,y_3)=b_2\psi_2(x_3,y_3)-\psi'_2(b_3x_3,b_3y_3)&
(x_3,y_3\in\ka), \\
  &
  \delta^1(f_1)(x_3,y_3)=b_1\psi_1(x_3,0,y_3)-\psi'_1\big(b_3x_3,f_2(y_3),b_3y_3\big)&
(x_3,y_3\in\ka).
\end{array}
    \end{equation}
Manifestly tori fixing the point $O$ correspond under $\Phi_1$; in
particular we have
$\Phi_1\big(T^{(e_2\!,\,e_3,\,h_2\!,\,h_3)}_\beta\big)=T^{\big(e'_2\!,\,e'_3,h'_2\!,\,h'_3\big)}_{\beta'}$
since tori are conjugated under $G_u$. This means
$$
(u_1, u_2, u_3)^{\Phi_1(\tau)}=\big(a_\tau^{\varepsilon e_1'}u_1,
a_\tau^{\varepsilon e_2'}u_2, a_\tau^{\varepsilon e'_3}u_3\big),
$$
 with $\varepsilon =\pm 1$. The identity
$\Phi_1\big((0,0,u_3)^\tau\big)=\big(\Phi_1(0, 0,
u_3)\big)^{\Phi_1(\tau)}$ and the first part of (\ref{eq2}) yield
$\varepsilon=1$, $e_3=e_3'$ and $f_j\big(a_\tau^{e_3}
u_3\big)=a_\tau^{e_j'} f_j(u_3)$, $j=1,2$, whereas
$\Phi_1\big((u_1,u_2,0)^\tau\big)=\big(\Phi_1(u_1,u_2,
0)\big)^{\Phi_1(\tau)}$ and (\ref{cinque}) give $e_1=e_1'$ and
$e_2=e_2'$. So the polynomials $f_j$ must be monomials and
consequently, in case $f_j\ne 0$,
\begin{equation}\label{DIECI}
    \begin{array}{cc}
\hbox{$e_j=e_3\!\deg(f_j)$}&(j=1,2).
    \end{array}
\end{equation}
Therefore $f_j(\mathrm{T})=d_j\mathrm{T}^\frac{e_j}{e_3}$,
$d_j\in\ka$, $j=1,2$. Furthermore  imposing the condition
$\Phi_1(0, 0, u_3) \Phi_1(0, v_2, 0)=\Phi_1\big((0, 0,u_3) (0,
v_2, 0)\big)$ we obtain $b_2^{h_2'}b_3^{h_3'}u_3^{h_3'} v_2^{h_2'}
= b_1u_3^{h_3} v_2^{h_2} $, i.e. $(h_2', h_3')=(h_2, h_3)$ and
\begin{equation}\label{XI}
    b_1=b_2^{h_2}b_3^{h_3}.
\end{equation}
So the first step is achieved:

\begin{all}{\bf $\bf\!\!.$ Proposition:}\label{XVI}
Let ${\sf G}^{\left(e'_2\!,\,e'_3,h'_2\!,\,h'_3\right)}_{\beta'}$
and ${\sf G}^{(e_2,\,e_3,\,h_2,\,h_3)}_\beta$ isomorphic as
algebraic permutation groups. Then
$$\big(e'_2,e'_3,h'_2,h'_3\big)=(e_2,e_3,h_2,h_3).$$\hfill$\Box$
\end{all}

Theorem 4.6 in \cite{DMG} says that the first of (\ref{f_3})
occurs precisely if
\begin{equation}\label{diciotto}
\delta^1(f_2)=  \psi_2-\psi'_2=0.
\end{equation}
 Also the fact that $e_1=e_3 \deg(\beta)$ if $\beta$ is not
the zero polynomial confines matters to examine  the case where
$\car =p>0$, $\beta=0$ and either
$\beta'(x_3,y_3)=\sum_{i=1}^{p-1}\frac 1
p\big(^{p\,}_{\,i}\big)x_3^{i{p^r}}y_3^{(p-i)p^{r}}\!\!,$ or
$\beta'(x_3,y_3)=x_3^{p^r}y_3^{p^s}$: by (\ref{DIECI}) we have
$\deg\beta'=\frac{e_1}{e_3}=\deg f_1$ in case $d_1\ne 0$. Then the
second identity of (\ref{f_3}) turns into
\begin{equation}\label{X}
    \delta^1 (f_1) (x_3, y_3)=-\beta'(b_3x_3, b_3y_3)- b_3^{h_3} x_3^{h_3}
f_2(y_3)^{h_2}
\end{equation}
and again Theorem 4.6 in \cite{DMG} excludes  the possibility that
$f_2$ is the zero polynomial. Then $f_2$ is an additive monomial
by (\ref{diciotto}) and (\ref{DIECI}) forces $e_2$ to be a
$p$-power. Thus, in view of the main theorem, both $e_1$ and
$\deg\beta'$, are the sum of two $p$-powers. So just the following
two possibilities can occur: either $\beta'(x_3, y_3)=x_3^{p^r}
y_3^{p^{s}}$, or $\car =2$ and $\beta'(x_3, y_3)=x_3^{2^r}
y_3^{2^r}$. Thus the main theorem gives either
$e_2h_2+{e_3h_3}={e_3}\!\left(p^r+p^{s}\right)$, or
$e_2h_2+{e_3h_3}=e_3 2^{r+1}$, which means that   the pair of
$p$-powers $(h_2,h_3)$ is one of the following
\begin{equation}\label{XII}
    \begin{array}{rlc}
      1. & (h_2,
h_3)=\big(\frac{e_3}{e_2}p^{r},p^s\big);&\qquad\qquad\qquad
\qquad\qquad\qquad\qquad\qquad\qquad\\
2. & (h_2, h_3)= \big(\frac{e_3}{e_2}p^{s},p^r\big);\\
      3. & (h_2, h_3)=\big(\frac{e_3}{e_2}2^{r},2^r\big).
    \end{array}
\end{equation}
 As the right side of (\ref{X}) must be a co-boundary,
(\ref{XII}.1) gives, (\ref{XII}.2) and (\ref{XII}.3), lead
respectively to
\begin{equation*}\label{XIII}
    \begin{array}{rll}
     1. &
d_1=b_3^{p^r+p^{s}}=b_3^{p^s}\!\!d_2^{\frac{e_3}{e_2}p^{r}},&\text{hence}
\; d_2=b_3^\frac{e_2}{e_3};\qquad\qquad\qquad
\qquad\\
       2. & \hbox{$d_1=0\,$ and
$\;b_3^{p^r+p^{s}}=-b_3^{p^{r}}\!\!d_2^{\frac{e_3}{e_2}p^{s}}$},&\text{hence}
\; d_2=-b_3^\frac{e_2}{e_3};\\
      3. &
\hbox{$\;b_3^{2^{r+1}}=b_3^{2^{r}}\!\!d_2^{\frac{e_3}{e_2}2^{r}}$},&\text{hence}
\; d_2=b_3^\frac{e_2}{e_3}.
    \end{array}
\end{equation*}
Now it is straightforward calculation to verify that, for any
$b_1,b_2,d_3\in\ka$, the maps{\footnotesize
\begin{equation}\label{14}
    \begin{array}{rll}
      \!\!\!1. & \!\!\!\!\!\! \left\{\!\!\!%
\begin{array}{ll}
\!\!G^{\big(e_2,\,e_3,\frac{e_3}{e_2}p^{r},\,p^s\big)}_0\rightarrow
G^{\big(e_2,\,e_3,\frac{e_3}{e_2}p^{r},\,p^s\big)}_{x^{p^r}\!y^{p^{s}}
},
\\
\!\!(u_1,u_2,u_3,a)\mapsto\Big(b_2^{\frac{e_3}{e_2}p^{r}}b_3^{p^s}
u_1+(b_3u_3)^{p^r+p^s},\,b_2u_2+(b_3u_3)^{\frac{e_2}{e_3}},
b_3 u_3,\,a\Big)\!;\\
\end{array}%
\right.      \\
   \!\!\!   2. &\!\!\!\!\!\!  \left\{\!\!\!%
\begin{array}{ll}
\!\!G^{\big(e_2,\,e_3,\frac{e_3}{e_2}p^{s},\,p^r\big)}_0\rightarrow
G^{\big(e_2,\,e_3,\frac{e_3}{e_2}p^{s},\,p^r\big)}_{x^{p^r}\!y^{p^{s}}
},
\\
\!\!(u_1,u_2,u_3,a)\mapsto\Big(\!b_2^{\frac{e_3}{e_2}p^{s}}b_3^{p^r}u_1,\,b_2u_2-(b_3u_3)^\frac{e_2}{e_3},
\,b_3 u_3,\,a\Big)\!;\\
\end{array}%
\right.     \\
 \!\!\!     3. &\!\!\!\!\!\!  \left\{\!\!\!%
\begin{array}{ll}
\!\!G^{\big(e_2,\,e_3,\frac{e_3}{e_2}2^{r},\,2^r\big)}_0\rightarrow
G^{\big(e_2,\,e_3,\frac{e_3}{e_2}2^{r},\,2^r\big)}_{x^{2^r}\!y^{2^{r}}
},
\\
\!\!(u_1,u_2,u_3,a)\mapsto\Big(\!b_2^{\frac{e_3}{e_2}2^{r}}b_3^{2^r}u_1+d_1u_3^{2^{r+1}},
\,b_2u_2+(b_3u_3)^\frac{e_2}{e_3},b_3 u_3,\, a\Big)\!;\\
\end{array}%
\right.
    \end{array}
\end{equation}}
$\!\!$are group isomorphisms in correspondence to the values
(\ref{XII}.$\,i$) of the pair of $p$-powers $(h_2,h_3)$.
Manifestly such isomorphisms supply isomorphisms for the
associated permutation groups. Summing up we have

\begin{all}{\bf $\bf\!\!.$ Theorem:}\label{17}
The integer parameters $e_2, e_3, h_2 ,h_3$ and the polynomial
$\beta$ determine uniquely the isomorphy class of the
$(2,2)$-imprimitive algebraic group $\sf G$, except the cases
where the pair $(h_2,h_3)$ takes one of the (integer) values {\rm
(\ref{XII}.$\,i$)}  which produces the corresponding isomorphisms
{\rm (\ref{14}.$\,i$)}.\B
\end{all}

\bk

{\footnotesize \noindent C. Bartolone, A. Di Bartolo, Dipartimento
di Matematica e Applicazioni, Universit\à di Palermo, Via
Archirafi 34, I-90123 Palermo, Italy  \\ E-mail: cg@math.unipa.it;
alfonso@math.unipa.it}\mk

{\footnotesize \noindent K. Strambach, Mathematisches Institut,
Universit\"at Erlangen-N\"urnberg, Bismarckstr. $1\frac 12$,
D-91054 Erlangen, Germany \\ E-mail: strambach@mi.uni-erlangen.de}

\end{document}